\documentstyle[11pt]{article}

\def\tabaddress#1{{\small\it\begin{tabular}[t]{c}#1 \\[1.2ex]\end{tabular}}}
\def\UPCMAT{\it Departamento de Matem\'atica Aplicada y Telem\'atica\\
   Universidad Polit\'ecnica de Catalu\~na\\
   Campus Nord, M\'odulo C-3\\
   C/ Gran Capit\'an s.n.\\
   E-08071 BARCELONA\\
   SPAIN}

\newtheorem{teor}{Theorem}
\newtheorem{prop}{Proposition}

\newtheorem{lem}{Lemma}
\newtheorem{definition}{Definition}

\def\beq{\begin{equation}}
\def\eeq{\end{equation}}
\def\bea{\begin{eqnarray}}
\def\eea{\end{eqnarray}}
\def\beann{\begin{eqnarray*}}
\def\eeann{\end{eqnarray*}}
\def\beasn{\begin{sneqnarray}}
\def\eeasn{\end{sneqnarray}}
\def\ben{\begin{enumerate}}
\def\een{\end{enumerate}}
\def\bit{\begin{itemize}}
\def\eit{\end{itemize}}
\def\dst{\(\displaystyle}
\def\proof{( {\sl Proof} )\quad}
\def\derpar#1#2{\frac{\partial{#1}}{\partial{#2}}}

\def\mapping#1{\mathrel{\mathop{\longrightarrow}\limits^{#1}}}

\def\qed{\ifvmode\removelastskip\fi
{\unskip\nobreak\hfil\penalty50\hbox{}\nobreak\hfil
\hbox{\vrule height1.2ex width1.2ex}\parfillskip=0pt
\finalhyphendemerits=0 \par\smallskip}}
\def\vf{{\aleph}}
\def\df{{\mit\Omega}}
\def\Lag{{\cal L}}

\def\d{{\rm d}}
\def\h{{\rm h}}

\def\Real{{\bf R}}
\def\Tan{{\rm T}}
\def\Lie{\mathop{\rm L}\nolimits}
\def\inn{\mathop{i}\nolimits}
\def\Cinfty{{\rm C}^\infty}
\title{Non-standard connections in classical mechanics}
\author{\sc Arturo Echeverr\'ia-Enr\'iquez,
   \\
{\sc Miguel C. Mu\~noz-Lecanda\thanks{{\bf e}-{\it mail}: MATMCML@MAT.UPC.ES}},
   \\
{\sc Narciso Rom\'an-Roy\thanks{{\bf e}-{\it mail}: MATNRR@MAT.UPC.ES}},
   \\
   \tabaddress{\UPCMAT}}
\begin{document}
\maketitle
\thispagestyle{empty}

\begin{abstract}
In the jet-bundle description of first-order classical field theories
there are some elements, such as the lagrangian energy and the
construction of the hamiltonian formalism, which require
the prior choice of a connection.
Bearing these facts in mind, we analyze
the situation in the jet-bundle description of time-dependent
classical mechanics. So we prove that this connection-dependence
also occurs in this case, although it is usually hidden
by the use of the ``natural'' connection given by the
trivial bundle structure of the phase spaces in consideration.
However, we also prove that this dependence is dynamically
irrelevant, except where the dynamical variation
of the energy is concerned. In addition, the relationship between
first integrals and connections is shown for a large enough
class of lagrangians.
\end{abstract}

\vfill \hfill
\vbox{\raggedleft AMS s.\,c.\,(1980): 58F05, 70H. PACS: 0320i}\null
\clearpage

\section{Introduction}

One of the most interesting lines of current research in mathematical physics
is the geometric formulation of first-order classical field theories,
which is achieved mainly using jet bundles
$J^1E\to E\to M$ and the geometrical structures
with which they are endowed
\cite{BSF-gcf}, \cite{CCI-91}, \cite{EM-91}, \cite{EMR-glft}, \cite{Gc-74},
\cite{GM-83}, \cite{GS-73}, \cite{GIMMSY-mm}, \cite{Gu-87},
\cite{He-dgcv},  \cite{LR-85}, \cite{Sa-87}.

Among all the relevant features observed when dealing with
these geometrical formulations, we wish to point out the following:
there are some dynamical elements of the theories
depending on the prior choice of a connection in the configuration bundle
$\pi\colon E\to M$: For instance:
\bit
\item
In the lagrangian formalism,
the definition of the {\sl density of lagrangian energy} and the
{\sl lagrangian energy function} \cite{EMR-glft}.
\item
The construction of the hamiltonian formalism of these theories,
in particular the {\sl hamiltonian function}
and the evolution equations
(see, for instance, \cite{CCI-91}).
\eit

It is also known that time-dependent mechanical systems
can be geometrically described using jet bundles.
Then $M=\Real$, $E=Q\times\Real$ and $J^1E=\Tan Q\times\Real$
(where $Q$ represents the configuration space of the system),
and consequently this can be considered as a particular situation
of field theories \cite{EMR-91}, \cite{Gc-74}. In this context,
since $Q\times\Real$ is a trivial bundle, it is canonically
endowed with a ``natural'' connection which is used
(when necessary) to define all the dynamical
and geometric elements of the theory; in particular,
the hamiltonian formalism and the energy lagrangian function.
A possible conclusion  of this approach is that,
for non-autonomous mechanical systems, the definition of
these dynamical elements does not depend on the choice of any connection.

The aim of this paper is to
make evident the influence of the choice of a connection
in the geometrical construction of some elements of the theory.
In order to achieve this, we will choose an arbitrary connection
in the bundle $Q\times\Real\to\Real$ and we will re-construct
the dynamics of the theory starting from this point.
In addition, we will study the relation among the descriptions
coming from different choices of a connection, and we will interpret
the results of this analysis from a dynamical point of view.
Connections on $Q\times\Real\to\Real$ can be understood as
special time-actions in the manifold $Q\times\Real$.
The standard one is by translations $(q,t)\stackrel{s}{\mapsto}(q,t+s)$.
Other actions correspond to non-standard connections.
In this sense changes in connections imply changes in the energy and
in the geometric elements of the theory.

The structure of the work is as follows:

First, we introduce the basic ideas about connections
in the bundle $Q\times\Real$ and the natural geometric elements
in the bundle $\Tan Q\times\Real\to Q\times\Real$.
Section 3 is devoted to presenting the lagrangian formalism
(using arbitrary connections)
and showing which of its elements are connection-depending.
We obtain results on the variation of the energy along the
motion of the system and on the relationship between first integrals
and connections. All these results generalize classical ones for
non-autonomous systems.
In the fourth section we construct the hamiltonian formalism for
non-autonomous systems, depending on the choice of an arbitrary connection
and, subsequently, we give a dynamical interpretation
of the results so obtained.
In section 5 a characterization of the lagrangian energy function,
based on variational principles, is given.
A final section is devoted to summarizing the conclusions
reached in the work.

All the manifolds and maps are $\Cinfty$.
Sum over repeated indices is understood.

\section{The 1-jet bundle of $\pi\colon Q\times\Real\to\Real$.
Geometric structures and connections}

The ideas in this section are known. We merely emphasize the differences
between the general situation and this particular one
in order to make the paper more readable and selfcontained.
See \cite{GHV-72} and \cite{Sa-89} as general references.

\subsection{Connections in $\pi\colon Q\times\Real\to\Real$}

Consider the bundle $\pi\colon Q\times\Real\to\Real$,
where $Q$ is an $n$-dimensional differentiable manifold
(the {\sl configuration space} of a physical system).
The $1$-jet bundle of sections of $\pi$ is
$\pi_1\colon \Tan Q\times\Real\to Q\times\Real$.
In fact, if $\phi =(\phi_Q,{\rm Id})$ is a local section of $\pi$
defined in a neighbourhood of $s\in\Real$, with
$\phi (s)=(q,s)$, then the $1$-jet equivalence class
of $\phi$ is determined by \dst (\Tan_s\phi_Q)\left(\frac{\d}{\d t}\right)\) ;
that is, an element of $\Tan_qQ$. Conversely, if $v\in\Tan_qQ$
and $s\in\Real$, there is an equivalence class of curves
$\phi_Q\colon(s-\epsilon ,s+\epsilon )\to Q$ with
\dst\Tan_s\phi_Q\left(\frac{\d}{\d t}\right)=v\) ;
so a $1$-jet equivalence class of local sections is defined and
$\phi =(\phi_Q,{\rm Id})$ is one of its representatives.

If $(U;q^\mu)$ is a local chart in $Q$,
then a local chart in $\Tan Q\times\Real$
is $(\pi_1^{-1}(U\times\Real );q^\mu ,t,v^\mu )$ where
$$
v^\mu ((q,s),v) =
\left(\frac{\d\phi^\mu(t)}{\d t}\right)_s
$$
$\phi\colon\Real\to Q\times\Real$ being
a representative of $((q,s),v)\in\Tan Q\times\Real$
and $\phi^\mu=q^\mu \circ \phi_Q$.
Of course, these coordinates $v^\mu$ are the physical velocities.

In order to introduce the connections, we must study
the tangent bundle of $Q\times\Real$. Observe that
there is a natural identification between
$\Tan (Q\times\Real )$ and
$\Tan Q\times\Tan\Real :=\pi_Q^*\Tan Q\oplus\pi^*\Tan\Real$ given by
$$
\begin{array}{lllll}
\psi & \colon & \Tan (Q\times\Real ) & \longrightarrow & \Tan Q\times\Tan\Real
\\
& & ((q,s),u)& \mapsto & ((q,s),\Tan_{(q,s)}\pi_Q(u)+\Tan_{(q,s)}\pi (u))
\end{array}
$$
where $\pi_Q\colon Q\times\Real\to Q$ is the natural projection.

But $\pi_Q^*\Tan Q$ is identified with ${\rm V}(\pi )$
(the vertical subbundle of $\Tan (Q\times\Real )$ with respect to $\pi$).
In fact, if $(q_o,s)\in Q\times\Real$ and $j_s\colon Q\to Q\times\Real$
is the $s$-injection defined by $j_s(q)=(q,s)$, then
${\rm V}_{(q_o,s)}(\pi )=\Tan_{q_o}j_s(\Tan_{q_o}Q)$.
So we have the natural spitting
$$
\Tan (Q\times\Real )={\rm V}(\pi )\oplus\pi^*\Tan\Real
$$
and $\pi^*\Tan\Real$ is called the {\sl horizontal subbundle}.
As a consequence, if $v\in\Tan_{(q,s)}(Q\times\Real )$,
we will write $v=v_Q+v_{\Real}$ in this splitting.

This natural splitting will be called the
{\sl standard connection} in the bundle
$\pi\colon Q\times\Real\to\Real$. The theory of connections describes
the possible splittings of this kind.

Taking into account this model and following \cite{Sa-89},
we prove the following statement:

\begin{prop}
The following elements on  $\pi\colon Q\times\Real\to\Real$
can be canonically constructed one from the other:
\ben
\item
A section of $\pi_1\colon\Tan Q\times\Real\to Q\times\Real$,
that is a mapping $\nabla\colon Q\times\Real\to\Tan Q\times\Real$
such that $\pi_1\circ\nabla ={\rm Id}_{Q\times\Real}$.
\item
A subbundle ${\rm H}(\nabla )$ of $\Tan (Q\times\Real )$
such that
\beq
\Tan (Q\times\Real )={\rm V}(\pi )\oplus{\rm H}(\nabla )
\label{split}
\eeq
\item
A semibasic $1$-form $\tilde\nabla$ on $Q\times\Real$
with values in $\Tan (Q\times\Real )$ (that is, an element of
$\Gamma (Q\times\Real ,\pi^*\Tan^*Q)\otimes\vf (Q\times\Real )$),
such that $\alpha\circ\tilde\nabla =\alpha$, for every
semibasic form $\alpha\in\df^1(Q\times\Real )$.
\een
\end{prop}
(We use the notation $\Gamma (A,B)$ for the set of sections
of the bundle $B\to A$).

\proof
\quad (1 $\Rightarrow$ 3)\quad
Let $\nabla$ be a section of $\pi_1$ and
$\phi =(\phi_Q,{\rm Id})$ a representative of
$\nabla (q,s)\in\Tan Q\times\Real$. We have the map
$$
\Tan_s\phi\circ\Tan_{(q,s)}\pi\colon\Tan_{(q,s)}(Q\times\Real )\longrightarrow
\Tan_{(q,s)}(Q\times\Real )
$$
Let $\tilde\nabla\colon\Tan (Q\times\Real )\longrightarrow\Tan (Q\times\Real )$
be the map defined in this way for every $(q,s)\in Q\times\Real$.
Obviously $\tilde\nabla$ vanishes on the $\pi$-vertical vector fields.

In order to prove the differentiability of $\tilde\nabla$,
take a local chart $(q^\mu ,t,v^\mu )$.
The local expression of $\nabla$ is
$\nabla (q,s) =(q,s,\gamma^\mu (q,s))$. Then, if
\dst X=f\derpar{}{t}+g^\mu\derpar{}{q^\mu}\in\vf (Q\times\Real )\) ,
we have
$$
\tilde\nabla ((q,s);X)=
f(q,s)\left(\derpar{}{t}+\gamma^\mu\derpar{}{q^\mu}\right)\Big\vert_{(q,s)}
$$
which is a differentiable vector field in $Q\times\Real$. So
$\tilde\nabla\in\Gamma (Q\times\Real ,\pi^*\Tan Q)\otimes\vf (Q\times\Real )$.

The local expression of $\tilde\nabla$ is
\beq
\tilde\nabla =\d t\otimes\left(\derpar{}{t}+\gamma^\mu\derpar{}{q^\mu}\right)
\label{le1}
\eeq
Then, if $\alpha\in\df^1(Q\times\Real )$ is a semibasic form, then the relation
$\alpha\circ\tilde\nabla =\alpha$ holds, as was desired.

\quad\quad (3 $\Rightarrow$ 1)\quad
Given $\tilde\nabla$, from its local expression \ref{le1},
we can construct local maps
$\nabla\colon U\times\Real\to\Tan U\times\Real$ given by
$\nabla (q,s):=(q,s,\gamma^\mu(q,s))$,
for every $(q,s)\in V=U\times\Real$. But,
since $\tilde\nabla$ is globally defined, its local representations agree
on the intersections of local charts.
In the same way, the local expressions of $\nabla$ can be glued
and we get a globally defined section of $\pi_1$.
Moreover, if $\tilde\nabla_1\not=\tilde\nabla_2$, then they
have different local representations on the neighbourhood of a point, and
hence they induce different sections $\nabla_1$ and $\nabla_2$.

\quad\quad (3 $\Rightarrow$ 2)\quad
Given $\tilde\nabla$, since it is a semibasic form we have
$\tilde\nabla\circ\tilde\nabla =\tilde\nabla$. Then we have the splitting
$$
\Tan (Q\times\Real )={\rm Ker}(\tilde\nabla )\oplus{\rm Im}(\tilde\nabla )
$$
But ${\rm Ker}(\tilde\nabla )={\rm V}(\pi )$, and we have
the splitting (\ref{split}) taking ${\rm H}(\nabla )={\rm Im}(\tilde\nabla )$.

\quad\quad (2 $\Rightarrow$ 3)\quad
Given the subbundle ${\rm H}(\nabla )$ and
the splitting (\ref{split}), we have the maps
$$
h_{\nabla}\colon\Tan (Q\times\Real )\to {\rm H}(\nabla )
\ , \
v_{\nabla}\colon\Tan (Q\times\Real )\to {\rm V}(\pi )
$$
called the {\sl horizontal} and {\sl vertical projections}
(we will use the same symbols $h_{\nabla}$ and $v_{\nabla}$
for the natural extensions of these maps to vector fields).
Then we can define a form $\tilde\nabla$ with values in
$\Tan (Q\times\Real )$ by
$$
\tilde\nabla ((q,s);u):=h_{\nabla}((q,s),u)
$$
Obviously, $\tilde\nabla$ is semibasic and
$\alpha\circ\tilde\nabla =\alpha$ for every semibasic form.
\qed

\begin{definition}
A {\rm connection} in the bundle $\pi\colon Q\times\Real\to\Real$
is one of the above mentioned equivalent elements.

Then ${\rm H}(\nabla )$ is called the {\rm horizontal subbundle} of
$\Tan (Q\times\Real )$ associated with the connection $\nabla$
and its sections {\rm horizontal vector fields}.
$\tilde\nabla$ is called the {\rm connection form}.
\end{definition}

In a local chart $(q^\mu ,t,v^\mu )$, the local expressions
of these elements are
$$
\tilde\nabla =
\d t\otimes\left(\derpar{}{t}+\gamma^\mu\derpar{}{q^\mu}\right)
\quad , \quad
{\rm H}(\nabla ) =
{\rm span}\left\{\derpar{}{t}+\gamma^\mu\derpar{}{q^\mu}\right\}
$$

For every vector field
\dst X \equiv X_{\Real}+X_Q \equiv f\derpar{}{t}+X_Q \in\vf (Q\times\Real )\)
(where $X_Q\in\Cinfty (Q\times\Real )\otimes\vf (Q)$),
the horizontal and vertical projections are given by
\beann
X_{h_{\nabla}}\equiv h_{\nabla}(X)&=&
f\left(\derpar{}{t}+\gamma^\mu\derpar{}{q^\mu}\right)
\\
X_{v_{\nabla}}\equiv v_{\nabla}(X)&=&
X_Q-f\gamma^\mu\derpar{}{q^\mu}
\eeann
and we have the splitting
$$
X=X_Q+f\derpar{}{t}=X-\tilde\nabla (X)+\tilde\nabla (X)=
X_{v_{\nabla}}+X_{h_{\nabla}}
$$

Moreover, we have the following result:

\begin{prop}
Every connection in the bundle $\pi\colon Q\times\Real\to\Real$
induces a canonical lifting
\beann
\vf (\Real ) & \longrightarrow & \vf (Q\times\Real )
\\
X &\mapsto& \bar X
\eeann
where $\bar X$ is a horizontal vector field.
\label{prol}
\end{prop}
\proof
Let $\nabla$ be a connection, $(q,s)\in Q\times\Real$
and and $\phi =(\phi_Q,{\rm Id})$ a representative of $\nabla (q,s)$.
If $X\in\vf (\Real )$, we define $\bar X$ by
$$
\tilde X(q,s):= \Tan_s\phi (X_s)
$$
 From the local expression of $\nabla$ we deduce that $\bar X$
is a $\Cinfty$-vector field and it is horizontal because
$$
\tilde\nabla ((q,s),\bar X(q,s))=
(\Tan_s\phi\circ\Tan_{(q,s)}\pi )(\bar X_{(q,s)})=
(\Tan_s\phi\circ\Tan_{(q,s)}\pi\circ\Tan_s\phi )(X_s)=
\Tan_s\phi (X_s)=\bar X_{(q,s)}
$$
\qed

But $\vf (\Real )$ has a global generator, \dst\frac{\d}{\d t}\) ;
so, given a connection $\nabla$, we can take its lifting:
$$
\overline{\frac{\d}{\d t}}\Big\vert_{(q,s)}=
\Tan_s\phi_Q\left(\frac{\d}{\d t}\Big\vert_s\right)+
\Tan_s{\rm Id}_q\left(\frac{\d}{\d t}\Big\vert_s\right)=
\Tan_s\phi_Q\left(\frac{\d}{\d t}\Big\vert_s\right)+
\derpar{}{t}\Big\vert_{(q,s)}
$$
Observe that \dst\Tan_s\phi_Q\left(\frac{\d}{\d t}\Big\vert_s\right)\)
is a section of the bundle $\pi^*_Q\Tan Q\to Q\times\Real$;
that is, a time-dependent vector field in $Q$.
Then we have:

\begin{prop}
A connection in the bundle $\pi\colon Q\times\Real\to\Real$
is equivalent to a time-dependent vector field in $Q$.
\label{vfc}
\end{prop}
\proof
Given a connection $\nabla$,
taking \dst Y=\overline{\frac{\d}{\d t}}-\derpar{}{t}\) ,
we have the desired vector field.

Conversely, given a time-dependent vector field
$Y\colon Q\times\Real\to\pi^*_Q\Tan Q$,
we have a connection defined by
$$
\begin{array}{ccccc}
\nabla&\colon&Q\times\Real&\longrightarrow&\Tan Q\times\Real
\\
& &(q,s)&\mapsto&(Y_{(q,s)},s)
\end{array}
$$
\qed

If $Y$ is the vector field induced by $\nabla$, then
the connection form is written as
$$
\tilde\nabla = \d t\otimes\left(\derpar{}{t}+Y\right)
$$

As we have pointed out above,
the trivial bundle $\pi\colon Q\times\Real\to \Real$
has a natural connection: the standard one $\nabla_0$,
with ${\rm H}(\nabla_0)=\pi^*\Tan\Real$. In this case
\dst\tilde\nabla_0=\d t\otimes\derpar{}{t}\) ,
the associated time-dependent vector field is $Y_0=0$
and the lifting induced by $\nabla_0$ is given by
\dst\overline{\frac{\d}{\d t}}\Big\vert_{(q,s)}=
\derpar{}{t}\Big\vert_{(q,s)}\) .

If we have another connection $\nabla$ with associated vector field $Y$,
this lifting is
$$
\overline{\frac{\d}{\d t}}\Big\vert_{(q,s)}=
\derpar{}{t}\Big\vert_{(q,s)}+Y_{(q,s)}
$$
then we can understand that the ``lines of time''
induced by this connection are the integral curves
of the vector field \dst\derpar{}{t}+Y:= \tilde Y\) ,
which is called the {\sl suspension of} $Y$ \cite{AM-78}.

 From now on, we will refer to {\sl non-standard connections}
in the bundle $\pi\colon Q\times\Real$
in those cases that differ from the standard one.

\subsection{Geometric elements}

In the lagrangian formalism the dynamics takes place in the manifold
$\Tan Q\times\Real$.
Then, in order to set it, we need to introduce some geometrical elements
of the bundle $\pi_1\Tan Q\times\Real\to Q\times\Real$
(see \cite{EMR-91}, \cite{Gc-74}, \cite{GS-73} and \cite{Sa-89}
for details).

{\bf The structural canonical 1-form} \cite{Gc-74}

We can define a 1-form $\vartheta$ in $\Tan Q\times\Real$,
with values in $\pi_1^*{\rm V}(\pi )$, in the following way:
$$
\vartheta (((q,s),u);X):=(\Tan _{((q,s),u)}\pi_1-
\Tan_{((q,s),u)}(\phi\circ\pi\circ\pi_1))(X_{(q,s)})
$$
where $\phi$ is a representative of $((q,s),u)\in\Tan Q\times\Real$.

$\vartheta$ is called the {\sl structure canonical form} of
$\Tan Q\times\Real$. Its local expression is
$$
\vartheta = (\d q^\mu-v^\mu\d t) \otimes\derpar{}{q^\mu}
$$

{\bf The vertical endomorphisms}.

Taking into account that $\pi_1\colon\Tan Q\times\Real\to Q\times\Real$
is a vector bundle and the fiber on $(q,s)\in Q\times\Real$
is $\Tan_qQ\times\{s\}$, there exists a canonical diffeomorphism
between the $\pi_1$-vertical subbundle and $\pi_1^*(\Tan Q\times\Real )$,
that is,
$$
{\rm V}(\pi_1)\simeq\pi_1^*(\Tan Q\times\Real )
\simeq\pi_1^*\pi_Q^*\Tan Q\simeq\pi_1^*{\rm V}(\pi )
$$
We denote by
$$
{\cal S}:\pi_1^*{\rm V}(\pi )\to {\rm V}(\pi_1)
$$
the realization of this isomorphism and we will use the same notation
${\cal S}$ for its action on the modules of sections of these bundles.

If $(q^\mu,t,v^\mu )$ is a local system of canonical coordinates
in $\Tan Q\times\Real$, by construction we have that
\dst{\cal S}\left(\derpar{}{q^\mu}\right)=\derpar{}{v^\mu}\) , and so
$$
{\cal S} = \xi^\mu\otimes \derpar{}{v^\mu}
$$
where $\{\xi^\mu\}$ is the dual basis of
\dst\left\{\derpar{}{q^\mu}\right\}\) in $\pi_1^*{\rm V}(\pi )$.

${\cal S}$ is an element of
$\Gamma (\Tan Q\times\Real,\pi_1^*{\rm V}^*(\pi ))
\otimes\Gamma (\Tan Q\times\Real,{\rm V}(\pi_1))$.
Taking into account that the structure form $\vartheta$ is an element of
${\mit\Omega}^1(\Tan Q\times\Real ,\pi_1^*{\rm V}(\pi_1)) =
{\mit\Omega}^1(\Tan Q\times\Real )\otimes
\Gamma (\Tan Q\times\Real,\pi_1^*{\rm V}(\pi ))$,
using the natural duality, by contracting ${\cal S}$ with $\vartheta$
we obtain an element
$$
{\cal V} := \inn ({\cal S})\vartheta \in
{\mit\Omega}^1(\Tan Q\times\Real ) \otimes
\Gamma (\Tan Q\times\Real,{\rm V}(\pi_1 ))
$$
whose local expression is
$$
{\cal V}=\left(\d q^\mu-v^\mu\d t\right)\otimes\derpar{}{v^\mu}
$$
Notice that ${\cal V}$ can be thought as a
$\Cinfty (\Tan Q\times\Real )$-module
morphism ${\cal V}\colon\vf (\Tan Q\times\Real)\to\vf (\Tan Q\times\Real)$
with image on the $\pi_1$-vertical vector fields.

${\cal S}$ and ${\cal V}$ are called the
{\sl vertical endomorphisms} of $\Tan Q\times\Real$.

Observe that ${\cal S}$ and ${\cal V}$ are sections of different bundles.
However, if we have a connection $\nabla$, ${\cal S}$ can be also understood
as an endomorphism on $\vf (\Tan Q\times\Real )$. In fact,
notice that the splitting
$\Tan (Q\times\Real )={\rm V}(\pi )\oplus{\rm H}(\nabla )$,
induced by the connection $\nabla$, has a dual counterpart
$\Tan^*(Q\times\Real )={\rm V}^*(\pi )\oplus{\rm H}^*(\nabla )$.
So, ${\rm V}^*(\pi )$ is identified by $\nabla$ with a subbundle of
$\Tan^*(Q\times\Real )$. With the same notation as above,
if $\{\xi^\mu\}$ is the dual basis of
\dst\left\{\derpar{}{q^\mu}\right\}\) in $\pi_1^*{\rm V}(\pi )$,
then this identification is given by
$$
\xi^\mu \mapsto \d q^\mu -\gamma^\mu\d t
$$
if
\dst\nabla=\d t\otimes\left(\derpar{}{t}+\gamma^\mu\derpar{}{q^\mu}\right)\) .

In the same way, we have the splitting
$$
\pi_1^*\Tan^*(Q\times\Real )=
\pi_1^*{\rm V}^*(\pi )\oplus\pi_1^*{\rm H}^*(\nabla )
$$
and the injection
$j_{\nabla}\colon \pi_1^*{\rm V}^*(\pi )\hookrightarrow
\pi_1^*\Tan^*(Q\times\Real )$.
But this last bundle is a subbundle of $\Tan^*(\Tan Q\times\Real )$
whose sections are the $\pi_1$-semibasic $1$-forms
in $\Tan Q\times\Real$. Hence, by means of this injection,
${\cal S}$ is an element of
$\df^1(\Tan Q\times\Real )\otimes\vf (\Tan Q\times\Real )$,
with values in the vertical vector fields,
which will be denoted by ${\cal S}^{\nabla}$.
Its local expression is
$$
{\cal S}^{\nabla} = (\d q^\mu -\gamma^\mu\d t)\otimes \derpar{}{v^\mu}
$$

Now, we are able to consider the difference
${\cal S}^{\nabla}-{\cal V}\in
\df^1(\Tan Q\times\Real)\otimes\vf (\Tan Q\times\Real )$
whose local expression is
$$
{\cal S}^{\nabla}-{\cal V} =(v^\mu-\gamma^\mu )\d t \otimes\derpar{}{t}
$$
which will be used henceforth in order to characterize the lagrangian energy.

\section{Lagrangian formalism. Connections and lagrangian energy functions}

A {\sl time-dependent lagrangian function} is a function
${\cal L} \in \Cinfty (\Tan Q\times\Real )$.
As it is known, we can construct the lagrangian
forms associated with ${\cal L}$ using
the geometrical structure of $\Tan Q\times\Real $.

\begin{definition}
The {\rm Poincar\'e-Cartan $1$ and $2$-forms}
associated with the lagrangian function ${\cal L}$
are the forms in $\Tan Q\times\Real $ defined by
$$
\Theta_{\Lag} := \d\Lag\circ{\cal V}+\Lag\d t
\quad , \quad
\Omega_{\Lag} := -\d\Theta_{\Lag}
$$
\end{definition}

For the coordinate expressions of the Poincar\'e-Cartan forms we obtain
\beann
\Theta_{\Lag} &=&
\derpar{{\cal L}}{v^\mu}(\d q^\mu-v^\mu\d t)+{\cal L}\d t =
\left({\cal L}-v^\mu\derpar{{\cal L}}{v^\mu}\right)\d t+
\derpar{{\cal L}}{v^\mu}\d q^\mu
\\
\Omega_{\Lag}&=&
-\d\left(\derpar{\Lag}{v^\mu}\right)\wedge\d q^\mu +
\d\left(\derpar{\Lag}{v^\mu}v^\mu -\Lag\right)\wedge\d t
\eeann
Observe that these elements do not depend on the connection.

As usual, we say that a lagrangian function
${\cal L}$ is {\sl regular \/} iff its associated form
$\Omega_{\Lag}$ has maximal rank, which
is equivalent to demanding that
\dst{\rm det}
\left(\frac{\partial^2{\cal L}}{\partial{v^\mu}\partial{v^\nu}}\right)\)
is different from zero at every point.

Assuming the regularity of $\Lag$,
the dynamics of the system is described by a vector field
$X_{\Lag}\in{\cal X}(\Tan Q\times\Real )$,
which is a {\sl Second Order Differential Equation} (SODE),
such that:
\beq
\inn(X_{\Lag})\Omega_{\Lag} = 0
\quad , \quad
\inn(X_{\Lag})\d t = 1
\label{dineq}
\eeq
As a consequence, the integral curves of $X_{\Lag}$ verify the
{\it Euler-Lagrange equations} and
$X_{\Lag}$ is obviously independent of the connection.

A very different picture arises when we try to define intrinsically the
{\sl lagrangian energy function}. If we consider the standard connection,
it can be obtained as follows: take the lifting of
\dst\frac{\d}{\d t}\) from $\Real$ to $\Tan Q\times\Real$
given by the connection, which will be denoted
as usual by \dst\derpar{}{t}\) , then
$$
{\rm E}_{\Lag}=-\inn\left(\derpar{}{t}\right)\Theta_{\Lag}
$$
being \dst{\rm E}_{\Lag} = \derpar{\Lag}{v^\mu} v^\mu -\Lag\)
its local expression.

Of course, this changes if we consider a non-standard connection
$\nabla$. Then it is more appropriate to define the lagrangian energy function
associated to $\nabla$ in the following way:

\begin{definition}
Let $\nabla$ be a connection in $\pi\colon Q\times\Real\to \Real$,
$\tilde Y\in\vf (Q\times\Real )$ the lifting of
\dst\frac{\d}{\d t}\) induced by $\nabla$ and
$j^1\tilde Y\in\vf (\Tan Q\times\Real )$ its canonical lifting.
The {\rm lagrangian energy function} associated with the lagrangian function
$\Lag$ and the connection $\nabla$ is
$$
{\rm E}^{\nabla}_{\Lag}=-\inn (j^1\tilde Y)\Theta_{\Lag}
$$
\label{energ}
\end{definition}

In a local system of coordinates we have
$$
j^1\tilde Y=\derpar{}{t}+\gamma^\mu\derpar{}{q^\mu}+
\left(\derpar{\gamma^\mu}{t}+v^\nu\derpar{\gamma^\mu}{q^\nu}\right)
\derpar{}{v^\mu}
$$
and
$$
{\rm E}_{\Lag}^\nabla = \derpar{\Lag}{v^\mu} (v^\mu-\gamma^\mu )-\Lag
$$
It is obvious from this expression that
the lagrangian energy is connection-depending.
In particular, when the standard connection is used, then
$\gamma^\mu =0$, which is the fact that ``hides''
the explicit dependence on the connection of the energy
in classical mechanics.

Notice that, since $\Theta_{\Lag}$ is $\pi_1$-semibasic, the energy function
does not, in fact, depend  on the extension of $\tilde Y$ from
$Q\times\Real$ to $\Tan Q\times\Real$.

Another characterization of the energy can be obtained
in the following way.
In the geometrical description of mechanics,
it is also usual to define the lagrangian energy using
the vertical endomorphism. In order to achieve this result,
consider the difference ${\cal S}^{\nabla}-{\cal V}\in
\df^1(\Tan Q\times\Real)\otimes\Gamma (\Tan Q\times\Real,{\rm V}(\pi_1 ))$
and its natural contraction with $\d\Lag$,
$\d\Lag\circ ({\cal S}^{\nabla}-{\cal V})\in\df^1(\Tan Q\times\Real)$,
whose local expression is
$$
\d\Lag\circ ({\cal S}-{\cal V}) =\derpar{\Lag}{v^\mu} (v^\mu-\gamma^\mu )\d t
$$
Then, it is easy to verify that
$$
{\rm E}^\nabla_{\Lag}=
\inn (j^1\tilde Y)(\d\Lag\circ ({\cal S}-{\cal V})-\Lag\d t)
$$
The $1$-form
$$
{\cal E}^\nabla_{\Lag}:=\d\Lag\circ ({\cal S}-{\cal V})-\Lag\d t
\in\df^1(\Tan Q\times\Real )
$$
is called the {\sl density of lagrangian energy}
associated with the lagrangian function $\Lag$ and the connection $\nabla$
and it is a $(\pi\circ\pi_1)$-semibasic form
whose local expressions is
$$
{\cal E}_{\Lag}^\nabla =
\left(\derpar{\Lag}{v^\mu} (v^\mu-\gamma^\mu )-\Lag\right)\d t
$$

At this point it is relevant to ask about the dynamical meaning
of this connection-depending energies.
Remember that in the standard case we have
\beq
X_{\Lag}({\rm E}_{\Lag})=-\derpar{\Lag}{t}
\label{coner}
\eeq
which is a relation between the rates of variation of
two functions with respect to different time actions:
the action through the flow of the dynamical vector field $X_{\Lag}$
and the one given by translations in time.
As we have seen above, every connection $\nabla$
induces a lifting $j^1\tilde Y$ of \dst\frac{\d}{\d t}\)~,
which is the infinitessimal generator of a new time action,
and then the equivalent result to (\ref{coner})
is given by the following:

\begin{teor}
Let $X_{\Lag}\in\vf (\Tan Q\times\Real )$ be the dynamical vector field
(solution of the equations (\ref{dineq})). Then
\beq
X_{\Lag}({\rm E}^{\nabla}_{\Lag})=-(j^1\tilde Y)\Lag
\label{nce}
\eeq
\end{teor}
\proof
Notice that, since $X_{\Lag}$ is a SODE,
${\cal V}(X_{\Lag})=0$; therefore
\beann
X_{\Lag}({\rm E}^{\nabla}_{\Lag})&=&
\inn (X_{\Lag})\d{\rm E}^{\nabla}_{\Lag}=
-\inn (X_{\Lag})\d\inn(j^1\tilde Y)\Theta_{\Lag}=
\inn (X_{\Lag})\inn(j^1\tilde Y)\d\Theta_{\Lag}-
\inn (X_{\Lag})\Lie (j^1\tilde Y)\Theta_{\Lag}
\\ &=&
\inn(j^1\tilde Y)\inn (X_{\Lag})\Omega_{\Lag}+
\inn ([j^1\tilde Y,X_{\Lag}])\Theta_{\Lag}-
\Lie (j^1\tilde Y)\inn (X_{\Lag})\Theta_{\Lag}
\\ &=&
-\Lie (j^1\tilde Y)\inn (X_{\Lag})(\d\Lag\circ{\cal V}+\Lag\d t)=
-\Lie (j^1\tilde Y)\inn ({\cal V}(X_{\Lag}))\d\Lag
-\Lie (j^1\tilde Y)(\Lag\inn(X_{\Lag})\d t)
\\ &=&
-(j^1\tilde Y)\Lag
\eeann
where we have taken into account that
$[j^1\tilde Y,X_{\Lag}]$ is $\pi_1$-vertical
(because $X_{\Lag}$ is a SODE and $j^1\tilde Y$
is a canonical lifting \cite{Cr-83})
and $\Theta_{\Lag}$ is a $\pi_1$-semibasic form,
so $\inn ([j^1\tilde Y,X_{\Lag}])\Theta_{\Lag} =0$.
\qed

Furthermore, there is a relevant relationship between the
first integrals of the dynamical vector field and the lagrangian
energy functions which is given in the next statement:

\begin{prop}
If $\Lag$ is a lagrangian function such that its associated
Legendre transformation
(see next section) is different from zero at every point,
then every first integral of the dynamical vector field $X_{\Lag}$
is the energy associated to some connection.
\end{prop}
\proof
Let $\nabla$ be a connection such that its associated vector field
$Y$ verifies $(j^1\tilde Y)\Lag =0$, then it follows that
${\rm E}_\Lag^{\nabla}$ is a conserved quantity.
Conversely, if $f\in\Cinfty (\Tan Q\times\Real )$ verifies
$X_\Lag (f)=0$, and we take the connection associated to any vector field
$Y$ such that $\inn (j^1\tilde Y)\Theta_\Lag =-f$.
But taking into account that $\Theta_{\Lag}$ is
$\pi_1$-semibasic and using the definition of the Legendre transformation
$F\Lag$, this equation
can be written as
\dst Y\circ F\Lag +\inn\left(\derpar{}{t}\right)\Theta_\Lag =-f\) ,
which obviously can be solved for $Y$ if $F\Lag$
is different from zero at every point.
Then we have that ${\rm E}^{\nabla}_\Lag =f$.
\qed

To conclude this section, we may point out
that the lagrangian formalism
of non-autonomous systems does not depend
on taking the standard or non-standard
connections, except in order to define the energy lagrangian function,
which does depend on this choice.
This is the same situation as for the lagrangian formalism
of field theories (see \cite{EMR-glft}).

\section{Hamiltonian formalism with non-standard connections}
\protect\label{dinin}

Next we are going to consider the hamiltonian formulation of the problem,
taking the bundles
$\Tan^*Q\times\Real\stackrel{\tau_1}{\longrightarrow}
Q\times\Real\stackrel{\pi}{\longrightarrow}\Real$.
The construction and local expressions of the
canonical elements are the usual ones
when the standard connection is used \cite{Gc-74}, \cite{EMR-91}, .
Now we study how this situation changes when a
non-standard connection is taken.

First of all, we may recall that the transition from the lagrangian formalism
to the hamiltonian one is done by means of the
{\sl Legendre transformation}, which is defined for non-autonomous mechanics
as in the time-independent case.
So ${\rm F}\Lag\colon\Tan Q\times\Real\to \Tan^*Q\times\Real$
is the map ${\rm F}\Lag := ({\rm F}\Lag_t,Id_{\Real})$; that is,
for every $(x,t) \in \Tan Q\times\Real$,
$$
{\rm F}\Lag(x,t) := ({\rm F}\Lag_t(x),t)
$$
where ${\rm F}\Lag_t\colon\Tan Q\to\Tan^*Q$
is the usual fiber derivative of the
restriction of ${\cal L}$ obtained by considering its value on every fixed $t$,
that is, ${\rm F}\Lag_t(x) := {\rm F}\Lag(x,t) \vert_{t=ctn}$.
In coordinates we have
$$
{\rm F}\Lag^*q^\mu = q^\mu \quad , \quad
{\rm F}\Lag^*p_\mu = \derpar{\Lag}{v^\mu}
\quad , \quad {\rm F}\Lag^*t = t
$$
Then, the lagrangian function ${\cal L}$ is said to be {\sl regular} iff
${\rm F}\Lag$ is a local diffeomorphism and it is called {\sl hyperregular}
iff ${\rm F}\Lag$ is a global diffeomorphism.

Now, let $\nabla$ be a connection and
$\tilde\nabla\in\df^1(\Real )\otimes\vf (Q\times\Real )$ its
associated connection form. We can associate to it a canonical form
defined as follows:

\begin{definition}
Consider $(\alpha ,s)\in\Tan^*Q\times\Real$ and
$u\in\Tan_{(\alpha ,s)}(\Tan^*Q\times\Real )$, then let
$\theta_{\nabla}\in\df^1(\Tan^*Q\times\Real )$ be defined by
$$
\theta_{\nabla}((\alpha ,s);u):=\alpha (v_{\nabla}(\Tan\tau_1(u)))
$$
\ben
\item
The form $\theta_{\nabla}$ is called the
{\rm Liouville $1$-form} associated with the connection $\nabla$.
\item
The form $\omega_{\nabla}:=-\d\theta_{\nabla}$ is called the
{\rm Liouville $2$-form} associated with the connection $\nabla$.
\een
\end{definition}

If $t$ denotes the usual global coordinate in $\Real$,
$(q^\mu )$ is a local system of coordinates in $Q$ and
$(q^\mu ,p_\mu )$ is the canonical system of coordinates in $\Tan^*Q$,
the coordinate expressions of these elements are the following:
\beann
\tilde\nabla&=&\d t\otimes\left(\derpar{}{t}+\gamma^\mu\derpar{}{q^\mu}\right)
\equiv\d t\otimes\left(\derpar{}{t}+ Y\right)
\\
\theta_{\nabla}\left( (\alpha ,s);\derpar{}{q^\mu}\right) &=&
\alpha\left( v\left(\derpar{}{q^\mu}\right)\right)=
\alpha\left(\derpar{}{q^\mu}\right) =p_\mu (\alpha )
\\
\theta_{\nabla}\left( (\alpha ,s);\derpar{}{p_\mu}\right) &=&
\alpha (v(0))=0
\\
\theta_{\nabla}\left( (\alpha ,s);\derpar{}{t}\right) &=&
\alpha\left( v\left(\derpar{}{t}\right)\right)=
\alpha (- Y)=\alpha\left(-\gamma^\mu\derpar{}{q^\mu}\right) =
-\gamma^\mu (q,s)p_\mu (\alpha )
\eeann
(where $(q,s)=\tau_1(\alpha ,s)$), therefore
\beann
\theta_{\nabla} &=&
p_\mu\d q^\mu-\gamma^\mu p_\mu\d t
\\
\omega_{\nabla} &=&
\d q^\mu\wedge\d p_\mu+\gamma^\mu\d p_\mu\wedge\d t+
p_\mu\derpar{\gamma^\mu}{q^\nu}\d q^\nu\wedge\d t
\eeann

Using these coordinate expressions, it is easy to prove that:

\begin{prop}
${\rm F}\Lag^*\theta_{\nabla}=\Theta_{\Lag}+{\rm E}_{\Lag}^{\nabla}\d t$.
\label{pro5}
\end{prop}

If $\nabla$ is a connection in $\pi\colon Q\times\Real\to \Real$
and $Y\in\vf (Q\times\Real )$ its associated vector field,
let ${\cal Y}\in\vf (\Tan^*Q\times\Real )$ be any extension of $\tilde Y$ to
$\Tan^*Q\times\Real$; that is, any map
${\cal Y}\colon\Tan^*Q\times\Real\to\Tan (\Tan^*Q\times\Real )$
making the following diagram commutative
$$
\begin{array}{ccc}
\Tan (\Tan^*Q\times\Real )&\mapping{\Tan\tau_1}&\Tan (Q\times\Real )
\\
{\cal Y} \Big\uparrow \Big\downarrow& &\Big\downarrow \Big\uparrow \tilde Y
\\
\Tan^*Q\times\Real&\mapping{\tau_1}&Q\times\Real
\end{array}
$$
Consider now the projection $\tau_{\Tan^*Q}\colon\Tan^*Q\times\Real\to\Tan^*Q$.
Let $\theta\in\df^1(\Tan^*Q)$ and $\omega :=-\d\theta\in\df^2(\Tan^*Q)$
be the canonical forms in $\Tan^*Q$, and
$\theta_0=\tau_{\Tan^*Q}^*\theta$, $\omega_0=\tau_{\Tan^*Q}^*\omega$.
Then we can state:

\begin{prop}
The relation between the canonical forms in $\Tan^*Q\times\Real$
and the Liouville forms is given by
$$
\theta_{\nabla}=\theta_0-\eta \quad ,\quad
\omega_{\nabla}= \omega_0+\d\eta
$$
where $\eta\in\df^1(\Tan^*Q\times\Real )$ is the form defined as
\beq
\eta :=\theta_0({\cal Y})\d t\equiv H\d t
\label{funH}
\eeq
\label{lifor}
\end{prop}
\proof
Taking into account the coordinate expression of $\theta_{\nabla}$,
we must prove that the local expression of $\eta$ is just
\dst\eta \equiv \gamma^\mu p_\mu\d t\) or, what is equivalent,
for every $(\alpha ,s)\in\Tan^*Q\times\Real$ with
$(q,s)=\tau_1(\alpha ,s)$, that
$$
H(\alpha ,s):=\theta_0((\alpha ,s);{\cal Y})\equiv
\gamma^\mu (q,s)p_\mu (\alpha )
$$
Then, bearing in mind the diagram
$$
\begin{array}{ccc}
\Tan^*Q\times\Real&\mapping{\tau_1}&Q\times\Real
\\
\tau_{\Tan^*Q}\Big\downarrow\ & &\ \Big\downarrow \pi_Q
\\
\ \Tan^*Q&\mapping{\tau}&Q\
\end{array}
$$
we have
\beann
\theta_0((\alpha ,s);{\cal Y}) &=&
\theta (\alpha ;\Tan_{(\alpha ,s)}\tau_{\Tan^*Q}{\cal Y}(\alpha ,s))=
\alpha (\Tan_\alpha\tau\Tan_{(\alpha ,s)}\tau_{\Tan^*Q}{\cal Y}(\alpha ,s))=
\alpha (\Tan_{(\alpha ,s)}(\tau\circ\tau_{\Tan^*Q}){\cal Y}(\alpha ,s))
\\ &=&
\alpha (\Tan_{(\alpha ,s)}(\pi_Q\circ\tau_1 ){\cal Y}(\alpha ,s))=
\alpha (\Tan_{(q,s)}\pi_Q\Tan_{(\alpha ,s)}\tau_1{\cal Y}(\alpha ,s))=
\alpha (\Tan_{(q,s)}\pi_Q Y(q,s))
\\ &=&
\alpha (q; Y(q,s))=p_\mu (\alpha )\gamma^\mu (q,s)
\eeann
\qed

If we assume that $\Lag$ is a hyperregular lagrangian; that is,
${\rm F}\Lag$ is a global diffeomorphism; then
it is a simple matter of calculation to obtain the following alternative
characterization of the function $H$
$$
H =\theta_0({\rm F}\Lag_*(j^1\tilde Y))=
{\rm F}\Lag^{{-1}^*}(({\rm F}\Lag^*\theta_0)(j^1\tilde Y))
$$
since ${\rm F}\Lag_*(j^1\tilde Y)$ is an extension of $\tilde Y$ to
$\Tan^*Q\times\Real$.

In order to introduce the dynamics we need the
{\sl hamiltonian function}. Taking into account that the physical
observable corresponding to this function is the {\sl energy},
we have:

\begin{definition}
The {\rm hamiltonian function} associated with the connection $\nabla$
and the lagrangian $\Lag$
is the function $\h^{\nabla}\in\Cinfty (\Tan^*Q\times\Real )$
such that
$$
{\rm E}^{\nabla}_{\Lag}={\rm F}\Lag^*\h^{\nabla}
$$
(As one can observe, the existence of such a function is assured
(at least locally) if we assume that $\Lag$ is a hyperregular
(or regular) lagrangian.
In other cases, under certain hypothesis on $F\Lag$, it is possible to prove
that ${\rm E}^{\nabla}_{\Lag}$ is also ${\rm F}\Lag$-projectable).
\end{definition}

Now we define:

\begin{definition}
The {\rm Hamilton-Cartan ($1$ and $2$)-forms}
associated with the connection $\nabla$ and
the hamiltonian function $\h^{\nabla}$ are
$$
\Theta_{\nabla} := \theta_{\nabla}-\h^{\nabla}\d t
\quad ; \quad
\Omega_{\nabla} := -\d\Theta_{\nabla} = \omega_{\nabla}+\d\h^{\nabla}\wedge\d t
$$
\end{definition}

The equations of motion ({\sl Hamilton equations}) are
$$
\inn (X_{\nabla})\Omega_{\nabla} =0
\quad ,\quad
\inn (X_{\nabla})\d t=1
$$
where $X_{\nabla}\in\vf (\Tan^*Q\times\Real )$ is the
dynamical vector field associated with
the hamiltonian function $\h^{\nabla}$.

At this point we must consider the influence on the dynamics of the
choice of a non-standard connection
in the construction of the hamiltonian formalism.

In order to do this, we take the geometrical and dynamical
elements of the hamiltonian formalism when the standard connection $\nabla_0$
in $Q\times\Real$ is used. In this case the
{\sl standard hamiltonian function}
$\h\in\Cinfty (\Tan^*Q\times\Real )$ is defined by
$$
{\rm E}_{\Lag}={\rm F}\Lag^*\h
$$
Then, the relation between the hamiltonian functions $\h^{\nabla}$
and $\h$ is given by the following result:

\begin{prop}
Let $H$ be the function defined by (\ref{funH}). Then
$$
\h^{\nabla}=\h -H
$$
\label{relhf}
\end{prop}
\proof
It suffices to prove it using systems of natural coordinates.
Then, since $H=p_\mu\gamma^\mu$, the result holds immediately,
taking into account the coordinate expressions of the functions
${\rm E}^{\nabla}_{\Lag}$, ${\rm E}_{\Lag}$ and the
Legendre transformation.
\qed

Next we construct the Liouville forms associated with $\nabla_0$.
Taking into account the proposition \ref{lifor} and the fact that
for the standard connection is $ Y_0=0$, we obtain
straightforwardly that these forms are just $\theta_0$ and $\omega_0$.

Finally, the Hamilton-Cartan forms associated with $\nabla_0$
are defined as above, and so they are
$$
\Theta_0 := \theta_0-\h\d t
\quad , \quad
\Omega_0 := \omega_0+\d\h\wedge\d t
$$
Therefore, the equations of motion are now
$$
\inn (X_0)\Omega_0 =0
\quad ,\quad
\inn (X_0)\d t=1
$$
where $X_0\in\vf (\Tan^*Q\times\Real )$ is the
dynamical vector field associated with
the standard hamiltonian function $\h$.

Then we have the following result:

\begin{prop}
The Hamilton-Cartan forms associated with the standard and
non standard connections are the same:
$$
\Theta_0=\Theta_{\nabla }
\quad , \quad
\Omega_0=\Omega_{\nabla }
$$
As a consequence the dynamical vector fields are equal:
$X_{\nabla}=X_0$.
\end{prop}
\proof
Taking into account
the propositions \ref{lifor} and \ref{relhf} and the expression (\ref{funH})
we obtain that
$$
\Theta_{\nabla}:=\theta_{\nabla}-\h^{\nabla}\d t =
\theta_0-\eta-(\h-H)\d t=\theta_0-\h\d t =\Theta_0
$$
The remainder of the statement is a straightforward consequence of this result.
\qed

In a local system of canonical coordinates we have
$$
\begin{array}{ccccc}
\Theta_{\nabla} &=& p_\mu\d q^\mu -\h\d t&=& \Theta_0
\\
\Omega_{\nabla} &=& \d q^\mu\wedge\d p_\mu +\d\h\wedge\d t &=& \Omega_0
\\
X_{\nabla} &=&
\derpar{}{t}+\derpar{\h}{p_\mu}\derpar{}{q^\mu}
-\derpar{\h}{q^\mu}\derpar{}{p_\mu} &=& X_0
\end{array}
$$

Hence we can conclude that, from the dynamical point of view
in the hamiltonian formalism,
no difference exists when we describe the system using the standard or
non-standard connections.

As a final remark, one can easily check that
$\Theta_{\Lag}=F{\Lag}^*\Theta_{\nabla}$. Therefore
$\Omega_{\Lag}=F{\Lag}^*\Omega_{\nabla}$ and
$X_{\nabla}=F{\Lag}_*X_{\Lag}$.

\section{On the characterization of the energy
by means of variational principles}

In definition \ref{energ} we have defined the lagrangian energy
function, justifying it geometrically. However,
it can be considered as an ``ad hoc'' definition.
Next we give another characterization of the energy
which is based on variational principles and justifies
the definition we have given.

First of all we need the following lemma:

\begin{lem}
Let $\beta\in\df^1(\Tan Q\times\Real )$ and $f\in\Cinfty (\Tan Q\times\Real )$.
For every differentiable curve
$\sigma\colon [a,b]\subset\Real\to Q$,
the following conditions are equivalent:
\ben
\item
$\bar\sigma^*(f\d t)=\bar\sigma^*\beta$.
\item
\dst\int_{\bar\sigma}f\d t =\int_{\bar\sigma}\beta\) .
\een
(where $\bar\sigma\colon [a,b]\subset\Real\to\Tan Q\times\Real$
denotes the canonical lifting of $\sigma$).
\label{lema}
\end{lem}
\proof
Trivially $1\ \Rightarrow\ 2$.

Conversely, if we suppose $1$ is not true, then there exists one curve
$\sigma\colon [a,b]\subset\Real\to Q$ with
$\bar\sigma^*(f\d t-\beta )\not= 0$
and hence there is $s\in [a,b]$ and a closed neighbourhood
$V$ of $s$ in $[a,b]$ such that, taking
$\gamma\colon V\to Q$ with $\gamma =\sigma\vert_V$, then
$$
\int_{\bar\gamma}(f\d t -\beta )\not=0
$$
so $2$ is false.
\qed

Next we state:

\begin{prop}
The lagrangian energy function introduced in definition $\ref{energ}$
is the unique function in $\Tan Q\times\Real$ verifying the condition
$$
\bar\sigma^*({\rm E}^{\nabla}_{\Lag}\d t)=
\bar\sigma^*({\rm F}\Lag^*\theta_{\nabla}-\Lag\d t)
$$
for every curve $\sigma\colon [a,b]\subset\Real\to Q$.
\label{condi}
\end{prop}
\proof
(Uniqueness):\quad Let $f$ and $g$ be two functions verifying this condition.
Obviously $\bar\sigma^*((f-g)\d t)=0$, but
$0=\bar\sigma^*((f-g)\d t)=(f-g)(\bar\sigma (t))\d t$,
for every $t\in [a,b]$. Hence, $(f-g)(\bar\sigma (t))=0$,
and this implies $f-g=0$, because every point in $\Tan Q\times\Real$
is in the image of some curve $\bar\sigma$.

\quad\quad (Existence):\quad
 From proposition \ref{pro5} we obtain
\beann
\bar\sigma^*({\rm F}\Lag^*\theta_{\nabla}-\Lag\d t)&=&
\bar\sigma^*(\Theta_{\Lag}+{\rm E}_{\Lag}^{\nabla}\d t+\Lag\d t)
\\ &=&
\bar\sigma^*(\d\Lag\circ{\cal V}-\Lag\d t
+{\rm E}_{\Lag}^{\nabla}\d t+\Lag\d t)=
\bar\sigma^*({\rm E}_{\Lag}^{\nabla}\d t)
\eeann
since $\bar\sigma^*(\d\Lag\circ{\cal V})=0$.
So, the energy function introduced in definition \ref{energ} satisfies
this condition.
\qed

 From the lemma \ref{lema} we obtain, for every curve
$\sigma\colon [a,b]\to Q$,
$$
\int_{\bar\sigma}{\rm E}_{\Lag}^{\nabla}\d t=
\int_{\bar\sigma}({\rm F}\Lag^*\theta_{\nabla}-\Lag\d t)
$$
therefore
\beann
\int_{\bar\sigma}\Lag\d t &=&
\int_{\bar\sigma}({\rm F}\Lag^*\theta_{\nabla}-{\rm E}_{\Lag}^{\nabla}\d t)=
\int_{\bar\sigma}{\rm F}\Lag^*(\theta_{\nabla}-\h^{\nabla}\d t)
\\ &=&
\int_{\bar\sigma}{\rm F}\Lag^*(\theta_0-\h\d t)=
\int_{\bar\sigma}{\rm F}\Lag^*\Theta_0=
\int_{{\rm F}\Lag\circ\bar\sigma}\Theta_0
\eeann
and this equality shows the equivalence between the
{\sl Hamilton principle of minimal action}
(of the lagrangian formalism)  and the
{\sl Hamilton-Jacobi principle} (of the hamiltonian formalism).
Therefore, taking into account the proposition \ref{condi},
we have to conclude that the energy (as it is defined in definition
\ref{energ}) is the only function that realizes the equivalence
between both variational principles.
This fact justifies the definition given above.

\section{Conclusions}

The geometric description of non-autonomous mechanics is usually given
using the natural connection induced by the trivial bundle structure of
the phase spaces $\Tan Q\times\Real$ (for the lagrangian formalism)
and $\Tan^*Q\times\Real$ (for the hamiltonian formalism).
We have reformulated both formalisms starting from
the choice of an arbitrary connection. As a consequence of this analysis,
we have shown that the geometric construction of
some elements of the theory depends on this choice; namely,
the lagrangian energy function of the lagrangian formalism
and the hamiltonian formalism itself; in particular its
geometrical structures (Liouville forms) and the
hamiltonian function. This fact is hidden in
the usual geometric descriptions because the natural connection
associated to the trivial bundle structures is flat.

We have generalized the geometric definition
of the lagrangian energy function in order to take into account the
use of non-standard connections.
This definition characterizes the lagrangian energy as the
only function that realizes the equivalence between
the Hamilton variational principle and the Hamilton-Jacobi principle.

The next step was to investigate
the dynamical relevance of the choice of different connections.
We have proved that the dynamics is insensitive to this choice
because the Poincar\'e-Cartan forms of the lagrangian formalism and the
Hamilton-Cartan forms of the hamiltonian one do not depend on
the chosen connection.

It is worth pointing out that, from the physical point of view,
the connection-dependence of the energy function
can be understood as follows:
for time-dependent systems, to take a non-standard connection
can be interpreted as a change of the time action on $Q\times\Real$
(which arises from taking ${\rm H}(\nabla )$ instead of
$\pi^*\Tan\Real$ as the horizontal subbundle).
Then, it is reasonable that the energy function changes in its turn,
since it can be considered as the conjugate function of ``time''
(as is made evident in some geometrical descriptions
of non-autonomous systems \cite{EMR-91}).
A significant result is the relationship between the dynamical variation
of the connection-dependent energy and the variation of the
lagrangian with respect to the time-action induced by the connection.

A further consequence is that,
for time-dependent dynamical systems described by a
given lagrangian function, we have proved that any
first integral of the dynamics is the energy function
for a suitable connection. Moreover, we have a means of obtaining different
lagrangian energy (or hamiltonian) functions which are dynamically equivalent.
This consists of taking different connections in order to construct these
functions.

As a final remark, we trust that these results will help to clarify some
aspects
of the geometrical description of classical field theories.

\subsection*{Acknowledgments}

We are grateful for the financial support
of the CICYT TAP94-0552-C03-02.
We thank Mr. Jeff Palmer for the correction of the English version
of the manuscript.


\begin{thebibliography}{cccc}

\itemsep 1pt plus 1pt

\bibitem{AM-78}
{\sc R. Abraham, J.E. Marsden},
{\sl Foundations of Mechanics\/} (2nd ed.),
     Addison-Wesley, Reading Ma., 1978.

\bibitem{BSF-gcf}
{\sc E. Binz, J. \'Sniatycki, H. Fisher},
{\sl The Geometry of Classical fields},
North Holland, Amsterdam, 1988.

\bibitem{CCI-91}
{\sc J.F. Cari\~nena, M. Crampin,  L.A. Ibort},
``On the multisymplectic formalism for first order field theories'',
\sl Diff Geom. Appl.
\bf 1 \rm (1991) 345-374.

\bibitem{CIL-ts}
{\sc J.F. Cari\~nena, L.A. Ibort, E.A. Lacomba},
``Time scaling as an infinitesimal canonical transformation'',
{\sl Cel. Mec.}
{\bf 42} (1988) 201-213.

\bibitem{Cr-83}
{\sc M. Crampin},
``Tangent bundle geometry for Lagrangian dynamics'',
{\sl J. Phys. A: Math. Gen.\/}
{\bf 16} (1983) 3755-3772.

\bibitem{EM-91}
{\sc A. Echeverr\'\i a-Enr\'\i quez, M.C. Mu\~noz-Lecanda},
``Variational calculus in several variables: a hamiltonian approach'',
{\sl Ann. Inst. Henri Poincar\'e} {\bf 56}(1) (1992) 27-47.

\bibitem{EMR-91}
{\sc A. Echeverr\'\i a-Enr\'\i quez, M.C. Mu\~noz-Lecanda, N. Rom\'an-Roy},
``Geometrical setting of time-dependent regular systems.
Alternative models'',
{\sl Rev. Math. Phys.} {\bf 3}(3) (1991) 301-330.

\bibitem{EMR-glft}
{\sc A. Echeverr\'\i a-Enr\'\i quez, M.C. Mu\~noz-Lecanda,
N. Rom\'an-Roy},
``Geometry of Lagrangian First-order Classical Field Theories'',
Preprint DMAT-UPC (1994).

\bibitem{Gc-74}
{\sc P.L. Garc\'ia}
``The Poincar\'e-Cartan invariant in the Calculus of Variations'',
{\sl Symp. Math.}
{\bf 14} (1974) 219-240.

\bibitem{GM-83}
{\sc P.L. Garc\'ia, J. Mu\~noz},
``On the geometrical structure of higher order variational calculus'',
{\sl Proc. IUTAM-ISIMM Symposyum on Modern Developments
in Analytical Mechanics}
M. Francaviglia, A. Lichnerowicz eds.,
{\it Atti della Academia della Scienze di Torino}
{\bf 117} (1983) 127-147.

\bibitem{GS-73}
{\sc H. Goldschmidt, S. Sternberg},
``The Hamilton-Cartan formalism in the calculus of variations'',
{\sl Ann. Inst. Fourier Grenoble} {\bf 23} (1),
(1973) 203-267.

\bibitem{GIMMSY-mm}
{\sc M.J. Gotay, J.Isenberg, J.E. Marsden, R. Montgomery,
     J. \'Sniatycki, P.B. Yasskin},
{\sl Momentum maps and classical relativistic fields},
GIMMSY, 1990.

\bibitem{GHV-72}
{\sc W. Greub, S. Halpering, S. Vanstone},
{\sl Connections, curvature and cohomology},
Acad. Press, New York, 1972.

\bibitem{Gu-87}
{\sc C G\"unther},
``The polysymplectic hamiltonian formalism in the field theory
and the calculus of variations I: the local case'',
{\sl J. Diff. Geom.} {\bf 25} (1987) 23-53.

\bibitem{He-dgcv}
{\sc R. Hermann},
{\sl Differential Geometry and the Calculus of Variations},
Academic Press, $2^{nd}$ ed., Math. Sci. Press, Brookline Ma., 1968.

\bibitem{LR-85}
{\sc M. de Le\'on, P.R. Rodrigues},
{\sl Generalized Classical Mechanics and Field Theory},
North-Holland Math. Studies {\bf 112},
Elsevier, Amsterdam, 1985.

\bibitem{Ra-91}
{\sc M.F. Ra\~nada},
``Extended tangent bundle formalism for time-dependent
Lagrangian systems'',
{\sl J. Math. Phys.} {\bf 32} (2) (1991) 500-505.

\bibitem{Sa-87}
{\sc D.J. Saunders},
``The Cartan form in Lagrangian field theories'',
{\sl J. Phys. A: Math. Gen.}
{\bf 20} (1987) 333-349.

\bibitem{Sa-89}
{\sc D.J. Saunders},
{\sl The Geometry of Jet Bundles},
London Math. Soc. Lect. Notes Ser.
{\bf 142}, Cambridge, Univ. Press, 1989.

\end{thebibliography}
\end{document}